\documentclass[12pt]{article}
\usepackage[all,poly]{xy}
\xyoption{all}
\CompileMatrices
\usepackage{amsfonts}
\usepackage{amssymb}
\usepackage{amsmath}
\usepackage{epsfig}
\usepackage{amscd}
\usepackage{amsthm}
\usepackage{verbatim}
\begin{document}
\title{Representations of Surface Groups and Right-Angled Artin Groups in Higher-Rank}
\author{Stephen Wang}
\date{}
\maketitle

\newcommand{\R}{\mathbb R}
\newcommand{\Hyp}{\mathbb H}
\newcommand{\Z}{\mathbb Z}
\newcommand{\C}{{\mathcal {C}}}
\newcommand{\Cx}{{\mathbb {C}}}
\newcommand{\Q}{\mathbb Q}
\newcommand{\F}{{\bf {F}}}
\newcommand{\lat}{\Lambda}
\newcommand{\pf}{\noindent {\it Proof.}\hskip.2in}
\newcommand{\eps}{\epsilon}

\newcommand{\inv}{^{-1}}

\newtheorem{thm}{Theorem}
\newtheorem{conj}{Conjecture}
\newtheorem{lemma}{Lemma}
\newtheorem{defn}{Definition}
\newtheorem{prop}{Proposition}
\newtheorem{qst}{Question}

\begin{abstract}
We give very flexible, concrete constructions of discrete and faithful representations
of right-angled Artin groups into higher-rank Lie groups.  Using the geometry of the 
associated symmetric spaces and the combinatorics of the groups, we 
find a general criterion for when discrete and faithful representations exist, and 
show that the criterion is satisfied in particular cases.  There are direct applications towards
constructing representations of surface groups into higher-rank Lie groups, and, in particular,
into lattices in higher-rank Lie groups.
\end{abstract}

\section{Introduction}

The study of embeddings of surface groups $\Gamma = \pi_1(\Sigma_g)$, where $\Sigma_g$ is a closed surface
of genus $g$, in Lie groups
has a 
long history, and has been recently the subject of study from a variety of different viewpoints.  

The most basic example is when the Lie group $G = PSL(2, \R)$, in which case the faithful, discrete 
representations inside ${\rm Hom}(\Gamma,G)/G$ give Teichmuller space.  
Hitchin \cite{hitchin} expanded the notion of Teichmuller space to representations into 
$G=SL(n, \R)$, showing in that case that ${\rm Hom}(\Gamma,G)/G$ has a distinguished connected component, topologically
a ball, 
containing Teichmuller space.  He also calculated the number of connected components in the whole representation
space.   
Other examples
include
real projective structures on surfaces \cite{choigoldman},
K\"{a}hler geometry and bounded cohomology \cite{biw}, and 
Anosov flows and hyperconvex curves in projective space \cite{labour}.  

Looking at the representation variety is also a key component in studying vector bundles over surfaces and
surface bundles over surfaces.  In addition, when the images of the representations lie in a lattice $\Lambda$ in $G$, 
one obtains examples of essential surfaces inside a locally symmetric space, generalizing the 
study of closed geodesics to higher dimensional submanifolds.  
However, while much has been proven about general properties of ${\rm Hom}(\Gamma,G)/G$, such as the number
of connected components, or the
nature of certain types of embeddings, there has not been very much attention paid towards constructing
explicit examples of such objects, save in some special cases.  

This paper focuses on the general problem of constructing explicit examples of discrete, faithful 
representations of surface
groups into all sorts of Lie groups, or, even better, into lattices.  One can then hope to calculate numerical
invariants of those representations, such as Toledo invariants (when the associated symmetric space is Hermitian),
to 
determine which components of the representation variety these representations belong to.  
We apply techniques
and results of geometric group theory, combined with an understanding of the geometry of symmetric spaces,
to the problem.  

The main part of this work is the construction of  numerous new 
explicit examples of discrete, faithful representations of surface groups into semisimple Lie groups of 
higher rank.  We also construct new representations into lattices in semisimple Lie groups.  
The method developed  can also be adapted to representations of other groups of interest in
geometric group theory, such as the fundamental groups of closed hyperbolic 3- and 4-manifolds
(see \cite{crispwiest1}, \cite{janusz}), and
graph braid groups (\cite{crispwiest2}).

We examine the topic of representations of surface groups using the technique of right-angled
Artin groups.  Given a finite graph $H = (V,E)$, we define the \emph{right-angled Artin
group} $A(H)$ to be the group given by the presentation $\langle s_v | v\in V;
[s_v,s_w] = 1 \ {\rm{if}} \ vw \in E \rangle$.  Despite their simple presentations, these
groups exhibit a number of interesting properties, and have been the subject of an
increasing amount of study in recent years (see, e.g., \cite{charney}).  

The focus on right-angled Artin groups is useful in studying surface groups
because of the work of John Crisp 
and Bert Wiest \cite{crispwiest2}, who have shown that 
all fundamental groups of closed, oriented surfaces
embed in some right-angled Artin group.  In particular, the
fundamental group of any closed oriented surface of genus at least 2
embeds into the Artin group $A(C_5)$, where $C_n$ denotes the cycle
on $n$ vertices.  

This is one reason why we primarily focus on Artin groups $A(H)$ where $H$ is a cycle of at least 5 vertices.  
Another reason is that, since cycles are subgraphs of almost every graph, restrictions on the representations 
of the corresponding subgroups (or the lack thereof) will
translate to information about the full group.  

One can show by simple dimension count arguments that representations of $A(C_5)$ 
exist, but exhibiting faithfulness
and, particularly, discreteness requires a stronger argument.  
Our construction is explicit, and gives a whole class of discrete, faithful
representations.  In particular, the proof of discreteness exploits the geometry of the 
symmetric spaces associated to Lie groups.

The paper is in two main parts.  In the first, we prove that if there is an arrangement
of the maximal flats of the symmetric space $X = G/K$ that mimics the structure of the graph $H$ in a
certain sense, then 
we can find a faithful and discrete representation of $A(H)$ into $G$.  In the second, we primarily focus 
on the case where $H$ is a cycle, and show
that such arrangements of flats exist for certain symmetric spaces.

This process yields the following main results: 

\begin{thm}\label{firstmain}
There are infinitely many conjugacy classes of discrete, faithful representations of $A(C_5)$ 
into $SL(n,\R)$ for $n\ge 3$.  
\end{thm}

This method ought to be generalizable to arbitrary symmetric spaces, but it sometimes requires
Artin groups on  larger cycles than $C_5$.  For certain symmetric spaces, we will not be able to use our method 
to embed $A(C_5)$, but only Artin groups on cycles
of even length.  It is known that these Artin groups contain some, but not all, surface groups (\cite{droms}).  

\begin{thm}
There are infinitely many conjugacy classes of discrete, faithful representations 
of $A(C_6)$ into $SO(3,2)$.
\end{thm}

It is essential that the Lie groups are of $\R$-rank at least two.  

\begin{thm}
If $G$ is a semisimple Lie group of $\R$-rank 1, and $H$ is a non-complete connected
graph, then 
there is no discrete and faithful embedding of $A(H)$ into $G$.
\end{thm}

It is interesting that this method can also be used to construct 
representations of surface groups into lattices in semi-simple Lie groups.  
These constructions yield examples of essential surfaces inside locally symmetric spaces.  

\begin{thm}
There are infinitely many conjugacy classes of representations of $A(C_5)$ into 
$SL(5, \Z)$.  
\end{thm}

As a result of Theorem \ref{firstmain}, we obtain explicit representations of surface groups.  

\begin{thm}
For any closed surface group $\Gamma = \pi_1(\Sigma_g)$ with $g\ge 2$
and $n\ge 3$, 
there are infinitely many conjugacy classes of discrete, faithful representations 
of $\Gamma$ into $SL(n,\R)$.
\end{thm}

There is a little bit of work required to keep the infinitely many conjugacy classes,  
but otherwise this is an immediate corollary.  Certain discrete and faithful representations of 
surface groups into  $PSL(3,\R)$ correspond to convex projective structures on the surface, and
have been extensively studied (see  \cite{choigoldman},\cite{goldman2}).  For the most
part, however, the techniques used in previous work have not yielded algebraically explicit
examples.  Our construction gives a very flexible way of constructing these representations
explicitly.

All examples constructed via this method 
lie in the connected component of the trivial representation, so these representations
do not lie in Hitchin's Teichmuller component or have non-zero Toledo invariant (when the target
symmetric space is Hermitian).  These are, however, the first class of explicit non-trivial examples of
representations, and it is hoped that the geometric nature of the construction (and, in particular, the
ease in which it yields discreteness) will be helpful in future study of the representation variety.

\vskip.3in
\noindent
{\bf {Acknowledgements.}}  I wish to thank Benson Farb, under whose tutelage
at the University of Chicago this paper was conceived and (mostly) written.  He provided many
useful comments and suggestions, and his encouragement and enthusiasm were 
equally vital.  I would also like to thank Dave Morris and Alan Reid for 
their valuable suggestions, and also Kyle Ormsby for his help with the 
diagrams.  

\section{Preliminaries}

\subsection{Symmetric Spaces}

Throughout this paper, $G$ will denote a connected, semisimple Lie group with no compact or 
Euclidean factors, $K$ a maximal compact subgroup of $G$, and $X = G/K$ its associated symmetric
space of non-compact type.  

The {\emph{rank}} of $X$ is the largest integer $r$ such that there exists a totally geodesic
subspace of $X$ isometric to $\R^r$.  
A symmetric space of non-compact type has rank 1 if and only if its sectional curvatures are negative.  
A {\emph{maximal flat}} in a symmetric space $X$ is a totally geodesic submanifold $F$ which is isometric to
$\R^r$, where $r$ is the rank of $X$.  Note that, under the identification of $T_pX$ with $\mathfrak{p}$, 
$T_pF$ is a maximal abelian subalgebra $\mathfrak{a}\subset \mathfrak{p}$.  

Every geodesic of $X$ is contained in at least one maximal flat.  
A geodesic $\gamma$ is called {\emph{regular}} if there is exactly one maximal flat containing $\gamma$;  
if there is more than one maximal flats containing $\gamma$, it is {\emph{singular}}.  
Likewise, a point $z \in X(\infty)$ is called regular (resp. singular) if a geodesic $\gamma$ with
$\gamma(\infty) = z$ is regular (resp. singular).  

An isometry $\phi$ of $X$ having the property that there is some geodesic $\gamma$ such that
$\phi(\gamma(t)) = \gamma(t + t_0)$ for all $t$ is called {\emph{axial}}.  In this case, we call 
$\gamma$ an axis for $\phi$ and $t_0$ the translation distance.  

For these and other facts on symmetric spaces, see (for instance) \cite{ebbook}.  

\subsection{Artin Groups}

Given a finite graph $H = (V,E)$, we define the {\emph{ right-angled Artin
group}} $A(H)$ to be the group given by $\langle s_v | v\in V;
[s_v,s_w] = 1 \ {\rm{if}} \ vw \in E \rangle$.

Crisp and Wiest have shown that all but three surface
groups embed in some right-angled Artin group.  In particular, for closed orientable surfaces, we have
the following:

\begin{thm}[Crisp-Wiest \cite{crispwiest2}]
There is a faithful homomorphism from the fundamental group of any closed oriented surface $\Sigma_g$ 
into the Artin group $A(C_5)$, where $C_n$ denotes the cycle
on $n$ vertices.  
\end{thm}

Other Artin groups also contain certain surface groups as subgroups:

\begin{thm}[Droms-Servatius-Servatius \cite{droms}]
If $n\ge 5$, then $A(C_n)$ contains the surface group $\pi_1(\Sigma_g)$, where $g = 1 +(n-4)2^{n-3}$.  
\end{thm}

We will therefore concentrate on the question of finding discrete and
faithful representations of right-angled Artin groups, particularly $A(C_5)$, into semi-simple 
Lie groups.

It will be important that the Lie group be of rank greater than 1.  
There are examples of faithful embeddings of $A(C_5)$ into Lie groups
of rank 1, but they are far from being discrete.

\begin{thm}
If $G$ is a semisimple Lie group of $\R$-rank 1, and $H$ is a non-complete connected
graph, then 
there is no discrete and faithful embedding of $A(H)$ into $G$.
\end{thm}

Let $v$ and $w$ be two adjacent vertices of $H$, and $s_v$ and
$s_w$ the corresponding generators of the right-angled Artin group.  
Given a representation $\sigma: A(H) \rightarrow G$, we examine the
action of $\sigma(s_v)$ on the symmetric space $X=G/K$, which will be of
negative curvature.  

If $\sigma(s_v)$ fixes some point $p\in X$, then 
$\sigma(s_v)$ lies in a subgroup of $G$ which is compact.  Since $s_v$ has
infinite order in $A(H)$, this means that $\sigma$ is either non-faithful or has
non-discrete image.  

Thus $\sigma(s_v)$ fixes exactly one or two points in $X(\infty)$.  If it fixes two,
then there is some geodesic $\gamma: \R \rightarrow X$ and some $t_0 \in \R$ such that
$\sigma(s_v)(\gamma(t)) = \gamma(t + t_0)$ for all $t$.  This will be the unique geodesic
that is fixed by $\sigma(s_v)$, and thus $\gamma$ will be fixed by $\sigma(s_w)$ as well.
There will be some $t_1 \in \R$ such that
$\sigma(s_w)(\gamma(t)) = \gamma(t + t_1)$ for all $t$.  

But this means that, for any $\eps>0$, we can find $k, l \in \Z - \{0\}$ such that 
$\sigma(s_v^k s_w^l)$ translates $\gamma(0)$ a distance less than $\eps$ (possibly 0), so
either $\sigma$ is not faithful or it does not have a discrete image.  

We are left with the case where $\sigma(s_v)$ fixes exactly one point $x$ at infinity.  
The commutation relations then force every element of the image of $\sigma$ to fix $x$,
so $\sigma$ is a representation into a horospherical subgroup $N$ of $G$.  But $N$ is a 
nilpotent group (see \cite{eberlein2}) and is therefore amenable; since
$A(H)$ contains a non-abelian free subgroup, it cannot embed faithfully into $N$.  
$\Box$

\section{Discrete and Faithful Representations}

Given a connected graph $H$ with no triangle subgraphs, we will show that 
the right-angled Artin group $A(H)$ embeds into a
higher-rank Lie group $G$ if one can find a configuration of
geodesics and flats in the symmetric space $X=G/K$ that mimics the graph $H$.  
That is to say, if we view the singular geodesics throuigh some point $p_0 \in X$ as the 
vertices of a graph $H'$,  and the maximal flats containing $p_0$ as the edges of the graph, 
we wish to find an induced subgraph of $H'$ which is isomorphic to $H$.  We will use the 
terminology of graph theory and say that two geodesics through $p_0$ are {\emph{adjacent}} if they
lie in a common maximal flat.  

There will be a few technical conditions, so the exact statement is as follows:

\begin{thm}\label{mainthm}
Let $X = G/K$ be a symmetric space of non-compact type of rank at least 2, and let $p_0 \in X$.
Suppose that for each vertex $v$ of $H$, we can find a
geodesic $\gamma_v: \R \to X$ with $\gamma_v(0) = p_0$ with the following properties:
\begin{itemize}
\item
If $vw$ is an edge of $H$, then $\gamma_v$ and $\gamma_w$ are adjacent, there is exactly one maximal flat 
$F_{vw}$ containing both $\gamma_v$ and $\gamma_w$, and 
the set of singular geodesics $S_{vw}$ contained in the copy of $\R^2 \subset F_{vw}$
determined by $\gamma_v$ and $\gamma_w$ is finite.
\item
If $vw$ is not an edge of $H$, $\gamma_v$ and $\gamma_w$ are not adjacent.
\item
If $v_1v_2$ and $w_1w_2$ are disjoint edges in $H$, then no elements of $S_{v_1v_2}$ and $S_{w_1w_2}$ are
adjacent, with the exception of $\gamma_{v_i}$ and $\gamma_{w_j}$ if $v_iw_j$ is an edge
of $H$.  
\item
If $vw$ and $wx$ are two edges of $H$ sharing an endpoint, then no element of $S_{vw}-\{\gamma_w\}$ is
adjacent to any element of $S_{wx}-\{ \gamma_w\}$.  

\end{itemize}

Then
there  are infinitely many conjugacy classes of discrete, faithful representations $\sigma: A(H) \to G$.
\end{thm}

The proof of the theorem is a modification of the ``ping-pong" technique.  
The classic ``ping-pong" argument was first used 
to prove that two
elements $\alpha, \beta$ in $SL(2,\R)$ generated a free subgroup. 
It employed four open sets $U_{\alpha}, U_{\alpha\inv}, U_{\beta}, U_{\beta\inv} \subset \Hyp^2(\infty)$ and
showed that a word in $\alpha$ and $\beta$ sent a base point $z \in \Hyp^2(\infty)$ to the
open set corresponding to the first letter in the word, utilizing the fact that $\alpha$
and $\beta$ had disjoint repelling and attracting fixed sets.  

Our proof follows the same spirit, but there are complicating factors due to the
fact that $A(H)$ is not free and the more intricate geometry of higher rank symmetric spaces.  
For instance, our open sets will necessarily intersect one another, and we will be unable to 
use separate open sets for a generator and its inverse.

We begin the proof with a lemma.  

\begin{lemma}
Let $\mathfrak{g}  = \mathfrak{k} + \mathfrak{p}$ be the Cartan decomposition given by 
a point $p_0 \in X$, and let $A \in \mathfrak{p}$ be a non-zero vector.  
Let $F$ denote the union of all maximal flats containing the geodesic $\gamma (t) = \exp(tA)p$.  
Let $z_0 \in X(\infty)$;
then $\displaystyle{\lim_{t\rightarrow \infty} \exp(tA)z_0}$ will exist, and lie in $F(\infty)$.  
\end{lemma}

\pf
Let $z$ be any accumulation point of the sequence $\{\exp (tA) z_0 \}$, and 
let $y = \gamma( - \infty)$.  
Since $X$ is non-positively curved, the function $f(t) = \angle_{\gamma(t)} (y,z_0)$
will be non-increasing, and therefore there is some $\beta$ such that 
$\displaystyle{\lim_{t\rightarrow -\infty} f(t) = \beta}$.  

We claim that $\angle_{\gamma(t_0)} (y, z) = \beta$ for all $t_0$.  Indeed, 
for any $\eps > 0$, $\beta -\eps \le \angle_{\gamma(-t)} (y, z_0) \le \beta$
for $t\gg 0$, so, applying the isometry $\exp((t+t_0)A)$, we see that
$\beta -\eps \le \angle_{\gamma(t_0)} (y, \exp((t+t_0)A)z_0) \le \beta$ for $t\gg 0$.  
Hence $\angle_{\gamma (t_0)} (y, z)$ must lie in the interval $[\beta-\eps,\beta]$.
Since $\eps$ was arbitrary, $\angle_{\gamma (t_0)} (y, z) = \beta$, and the 
claim is proven.  

That the claim implies the lemma is precisely Claim E5 in \cite{eberlein}.  
$\Box$

Let $\{e_i\}$ denote the set of edges of $H$, and $S_i$ the corresponding
sets of singular geodesics defined in the statement of the theorem, and denote
the unique maximal flat containing the geodesics of $S_i$ be denoted by $F_i$.  
Let $S$ be the union of all $S_i$, and let $F(\gamma)$ denote the union of all
flats containing $\gamma$.

Two geodesics  $\gamma, \eta$ are adjacent if and only if $F(\gamma)(\infty)$
and $F(\eta)(\infty)$ are disjoint.  Therefore,  
for each geodesic $\gamma \in S$, we can find an 
open set $U_{\gamma} \subset X(\infty)$ containing $F(\gamma)(\infty)$ 
such that $\overline{U_{\gamma}} \cap \overline{U_{\eta}} \ne 
\emptyset$
iff there is some $i$ such that $\gamma, \eta \in S_i$.  
Let $U_i = \bigcap_{\gamma \in S_i} U_{\gamma}$; this will be an open set
containing $F_i(\infty)$.  

Let $z \in X(\infty)$ be a point not contained in the closure of any of the $U_{\gamma}$.

\begin{lemma}\label{axes}
It is possible to find real numbers $t_v$ such that the isometries $\phi_v = \exp(t_vA_v)$
have the following properties:

\begin{itemize}

\item
If $\gamma \in S$ is not adjacent to $\gamma_v$, then $\phi_v^k(U_{\gamma}) \subset U_{\gamma_v}$
for all $k \in \Z^{\times}$.  

\item
For any $k \in \Z^{\times}$, $\phi_v^k(z) \in U_{\gamma_v}$.  

\item
If $e_i = vw$ is an edge in $H$, and $\gamma \in S$ is adjacent to neither $\gamma_v$ nor $\gamma_w$,
then for any $k,l \in \Z^{\times}$, 
$\phi_v^k \circ \phi_w^l (U_{\gamma}) \subset \bigcup_{\eta \in S_i} U_{\eta}$.

\item
If $e_i = vw$ is an edge in $H$, then for any $k,l \in \Z^{\times}$, 
$\phi_v^k \circ \phi_w^l (z) \in \bigcup_{\eta \in S_i} U_{\eta}$.

\end{itemize}
\end{lemma}

The first two conditions are satisfied as long as the $t_v$ are large enough, 
by the lemma.  

To show that the second two are satisfiable, for a given edge $e_i = vw$, 
let $\mathfrak{h} \subset \mathfrak{p}$ be the subspace generated by $A_v$ and $A_w$, and let
$\tau: \mathfrak{h} \rightarrow X(\infty)$ be
given by $\tau(A) = \exp (A)(z')$, for some 
$z'\notin \bigcup_{\eta \in S_i} U_{\eta}$.  We will show that $\tau \inv (\bigcup_{\eta \in S_i} U_{\eta})$
has a bounded open set in $\mathfrak{h}$ as its complement.  This will suffice, for then we can simply pick $t_v$ and $t_w$ 
such that $\tau \inv (\bigcup_{\eta \in S_i} U_{\eta}) \subset \{ aA_v + bA_w \vert |a|>1 \ \mathrm{or} \ |b|>1 \}$.  

To do so, we simply apply the lemma to each geodesic $\exp (tA)(p_0)$, where $A\in \mathfrak{h}$ has length 1.  
All but finitely many will be regular geodesics, and for these there is some minimum $T_A$ such that $\exp(tA)(z') \subset
U_i$ whenever $|t|>T_A$.  If $\exp (tA)(p_0) = \eta \in S_i$, then we can still find some $T_A$ such that 
$\exp(tA)(z') \subset U_{\eta}$ when $|t|>T_A$.  Thus along any one-dimensional subspace of $\mathfrak{h}$,
the image of $\tau$ will eventually land in $\bigcup_{\eta \in S_i} U_{\eta}$.  The function $A\rightarrow T_A$
will be continuous, and thus there is some $T$ such that $\tau(A) \subset
\bigcup_{\eta \in S_i} U_{\eta}$  whenever $|A|>T$.
$\Box$

Once we have chosen our $\phi_v$, we then define the representation $\sigma: A(H) \rightarrow G$ by $\sigma(s_v) = \phi_v$.  
If $vw$ is an edge of $H$, we know that $\gamma_v$ and $\gamma_w$ are adjacent, and thus
$\phi_v$ and $\phi_w$ commute.  Therefore the relations of $A(H)$ are satisfied, and this
is a group homomorphism.  To show that the image of $\sigma$ is faithful and discrete, we 
will show that if $h$ is a non-identity element of $A(H)$, then 
$\sigma(h)(z) \in \bigcup_{\eta \in S_i} U_{\eta}$ for some $i$.

For each element $h \in A(H)$, let $\ell(h)$ denote the minimum integer $n$ such that
$h$ has a representation $h = s_{v_1}^{k_1}s_{v_2}^{k_2}\ldots s_{v_n}^{k_n}$
as a product of generators of $A(H)$ (where the $k_i$ can be either positive or negative
integers).  

By induction on $n$, we will show that any element $h$ has a representation as a word
in the generators $h = s_{v_1}^{k_1}s_{v_2}^{k_2}\ldots s_{v_n}^{k_n}$ where $n = \ell(h)$ and
either $\sigma(h)(z) \in U_{v_1}$ if $v_1v_2$ is not an edge of $H$, or $\sigma(h)(z) \in 
\bigcup_{\eta \in S_i} U_{\eta}$ if $v_1v_2 = e_i$ is an edge of $H$.  

When $\ell(h) = 1$, this is true by our choice of the $\phi_v$.  Now assume that
this is true for all $n<m$, and suppose that $\ell(h) = m$.  
Pick a representation $h = s_{v_1}^{k_1}s_{v_2}^{k_2}\ldots s_{v_m}^{k_m}$.  
By transposing commuting generators, it might be possible to switch $s_{v_1}^{k_1}$ to a
position farther right; suppose that the farthest right it can go is position $j$, in the
word $h=s_{w_1}^{l_1}s_{w_2}^{l_2} \ldots s_{w_m}^{l_m}$.  Thus we know that
$s_{w_i}$ commutes with $s_{w_j}$ for all $i<j$, and since there are no triangles in $H$,
$s_{w_i}$ does not commute with $s_{w_{i+1}}$ if $i<j-1$.

Also, by our induction hypothesis,
we can transpose commuting generators so that the word $h_j = s_{w_{j+1}}^{l_{j+1}} \ldots s_{w_m}^{l_m}$
has the property that  $\sigma(h_j)(z)$ is in $U_{w_{j+1}}$ if $w_{j+1}w_{j+2}$ is not an edge of $H$, or 
$\bigcup_{\eta \in S_i} U_{\eta}$ if $w_{j+1}w_{j+2}$ is the edge $e_i$ of $H$.  We know that $w_jw_{j+1}$
cannot be an edge of $H$ (else we could do another transposition to get the original left-most generator
$s_{v_1}$ farther right). 
Therefore, if $j=1$, we must have $\sigma(h)(z) \in U_{w_1}$, and we are done.

If $j\ge 2$,
we see that $\phi_{w_{j-1}}^{l_{j-1}}\circ \phi_{w_j}^{l_j}$ must take $\sigma(h_j)(z)$
into $\bigcup_{\eta \in S_i} U_{\eta}$, where $w_{j-1}w_j$ is the edge $e_i$ of $H$.  
If $j = 2$, we are then done.  If $j>2$, 
since $w_iw_{i+1}$ is
not an edge of $H$ for $i<j-1$, we can therefore conclude that $\sigma(h)(z) \in U_{w_1}$.  
This concludes the induction argument.

It remains to be shown that we can get infinitely many conjugacy classes.  In the proof
of Lemma \ref{axes}, we were free to choose the $t_v$, the translation length of the axial
isometries $\phi_v$, to be any numbers, as long as they were sufficiently large.  Given finitely
many discrete and faithful representations $\sigma_i$, there are only countably many translation lengths of axial
isometries in their images, so as long as we choose a $t_v$ distinct from all of them, we are
guaranteed a new conjugacy class.  

\section{Finding Geodesic Configurations}

Since the result of Crisp and Wiest holds for a right-angled Artin group on a 5-cycle, 
we are interested in finding arrangements of 5 singular geodesics that have the properties
required by the theorem in the previous chapter.

\newcommand{\diag}{ {\mathrm{diag}}}

\subsection{The Model Group}

\begin{thm}\label{sl3}
Let $X = SL(3, \R)/SO(3)$, and $p_0 \in X$.  
Then there are infinitely many conjugacy classes of sets of 5 geodesics through $p_0$
satisfying the conditions of Theorem \ref{mainthm} for the graph $H=C_5$.  
\end{thm}

Our goal is to find a different flat $F_{12}$ which intersects the original flat $F_{01}$
in one of the singular geodesics $\gamma_1$, choose a different singular geodesic $\gamma_2$ in $F_{12}$, find 
another flat $F_{23}$ which intersects $F_{12}$ in $\gamma_2$, and so on, until we are
able to find a flat $F_{40}$ which intersects $F_{01}$ in a geodesic $\gamma_0$, completing the cycle.  
To find successive flats, we will apply ``rotations" around the singular geodesics - that is,
isometries fixing the geodesic pointwise - to move one flat to the next (see diagram).  
We will see that the condition that this process cycle back around to the original flat 
essentially is a condition on the product of the chosen rotations.  

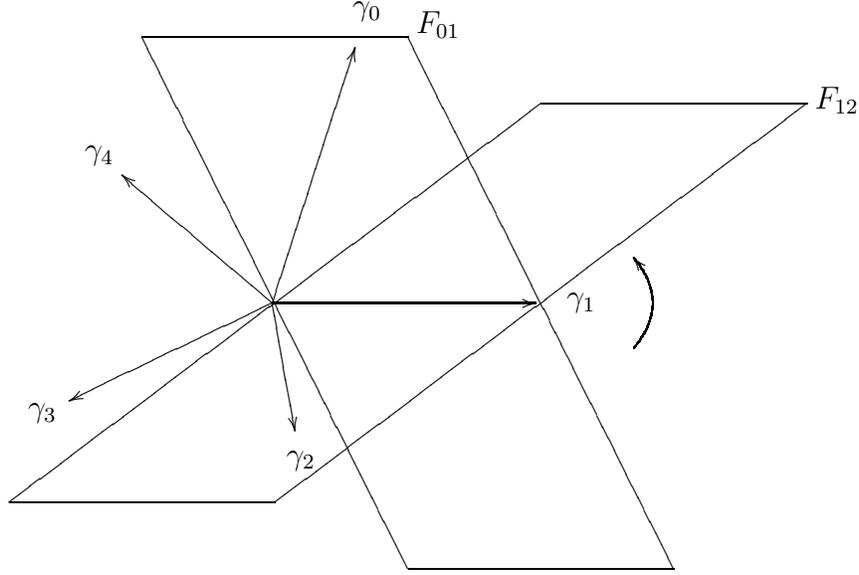
\begin{figure}[htbp]
  \begin{center}
    \[
    \begin{xy}/r1cm/:
    	(0,0),{\xypolygon4"F01"{~:{(2.5,0):(-1,2)::}}}
    	,(0,0),{\xypolygon4"F12"{~:{(2.5,0):(2,1.5)::}}}
    	,(-1.8,0);(-.7,3.4) **@{-}
    	?>*@{>} ?>*!/-5mm/{\gamma_0}
    	,(-1.8,0);(1.7,0) **@{-}
    	?>*@{>} ?>*!/-6mm/{\gamma_1}
    	,(-1.8,0);(-1.5,-1.7) **@{-}
    	?>*@{>} ?>*!/-4mm/{\gamma_2}
    	,(-1.8,0);(-4.5,-1.3) **@{-}
    	?>*@{>} ?>*!/-4mm/{\gamma_3}
    	,(-1.8,0);(-3.8,1.7) **@{-}
    	?>*@{>} ?>*!/-4mm/{\gamma_4}
    	,(.4,3.7)*{F_{01}}
    	,(5.7,2.7)*{F_{12}}
    	,(3,-.6);(3,.6) **\crv{(3.5,0)}
    	?>*@{>}
    \end{xy}
\]
    \caption{Flats and geodesics in a five-cycle configuration}
  \end{center}
\end{figure}

Without loss of generality, we may take $p_0$ to be the identity coset $K=SO(3)$.  

Consider the geodesics $\lambda_1(t) = \diag (e^{-2t}, e^t, e^t) p_0$, 
$\lambda_2(t) = \diag (e^t, e^{-2t}, e^{t}) p_0$, 
and $\lambda_3(t) = \diag (e^t, e^t, e^{-2t}) p_0$.  These 
are all contained in a unique common flat $F_{01}$.  
Let $A_i = \{ g\in K | g\lambda_i = \lambda_i \}$ for each $i$.  Note that
$A_1$ is the set of rotations in $SO(3)$ that fix the $x$-axis, $A_2$ fixes the
$y$-axis, and $A_3$ the $z$-axis.

Pick any non-identity elements $R_1 \in A_1$ and $R_2 \in A_2$.   

Let $Z$ be the vector $(0,0,1)$; consider $Y=R_2R_1Z$.  There is then an 
$R_3 \in A_3$ such that $R_3Y$ has second coordinate equal to 0, and thus
a unique $R_4 \in A_1$ such that $R_4R_3R_2R_1Z = Z$.  

Let $\gamma_0 = \lambda_2$, $\gamma_1 = \lambda_1$, 
$\gamma_2 = R_4\lambda_3$, $\gamma_3 = R_4R_3\lambda_2$ and 
$\gamma_4 = R_4R_3R_2\lambda_1$.  

Note that $\gamma_1$ and $\gamma_2$ both lie in the flat $F_{12} = R_4F_{01}$ 
since $R_4$ leaves $\gamma_1$ invariant.  Similarly, $\gamma_2$ and $\gamma_3$
both lie in the flat $F_{23} = R_4R_3F_{01}$, and $\gamma_3$ and $\gamma_4$ 
both lie in $F_{34} = R_4R_3R_2F_{01}$.  Since $R_4R_3R_2R_1 \in A_2$, 
both $\gamma_4$ and $\gamma_0$ are in $F_{40} = R_4R_3R_2R_1F_{01}$.  

Thus the five geodesics $\gamma_i$ satisfy the first condition of Theorem \ref{mainthm}.  
There are finitely many non-adjacency requirements; to verify that geodesics 
$\eta_1$ and $\eta_2$ are non-adjacent, we simply need to confirm that
axial isometries with $\eta_1$ and $\eta_2$ as axes do not commute.  

There are ten singular geodesics through $p_0$ in the five flats, namely
the five $\gamma_i$, plus $\eta_0 = \lambda_3$ in $F_{01}$, $\eta_1 = R_4\lambda_2$ in $F_{12}$, 
$\eta_2 = R_4R_3\lambda_1$ in $F_{23}$, $\eta_3 = R_4R_3R_2\lambda_3$ in $F_{34}$, and 
$\eta_4 = R_4R_3R_2R_1\lambda_3$ in $F_{40}$.  Let $T_1 = \diag ({{1}\over{4}}, 2, 2)$, 
$T_2 = \diag (2, {{1}\over{4}}, 2)$, and $T_3 = \diag (2, 2, {{1}\over{4}})$.  
Table \ref{sl3table} gives axial isometries for
each of the geodesics.

\begin{table}[h]\label{sl3table}
\begin{center}
\begin{tabular}{c|c}
	Geodesic & Isometry \\
        \hline
        $\gamma_0$ & $T_2$ \\
        $\gamma_1$ & $T_1$ \\
        $\gamma_2$ & $R_4 T_3 R_4\inv$ \\
        $\gamma_3$ & $R_4 R_3 T_2 R_3\inv R_4\inv$ \\
        $\gamma_4$ & $R_4 R_3 R_2 T_1 R_2\inv R_3\inv R_4\inv$ \\
	$\eta_0$ & $T_3$ \\
        $\eta_1$ & $R_4 T_2 R_4\inv$ \\
        $\eta_2$ & $R_4 R_3 T_1 R_3\inv R_4\inv$ \\
        $\eta_3$ & $R_4 R_3 R_2 T_3 R_2\inv R_3\inv R_4\inv$ \\
        $\eta_4$ & $R_4 R_3 R_2 R_1 T_3 R_1\inv R_2\inv R_3\inv R_4\inv$ 
\end{tabular}
\caption{Singular geodesics and axial isometries}
\end{center}
\end{table}

Pairs of geodesics that we need to be non-adjacent are $\gamma_i$ and $\gamma_j$ if $j-1\ne \pm 1$ (mod 5), 
$\eta_i$ and $\eta_j$ for any $i \ne j$, and $\eta_i$ and $\gamma_j$ if $j-i \ne 0, 1$ (mod 5).  

Recall that we have complete freedom to choose $R_1$ and $R_2$ in $A_1$ and $A_2$ respectively, but once 
they are chosen, $R_3$ and $R_4$ are determined.  Since commuting is an analytic condition,
if we see that these pairs of isometries do not commute for one particular choice of $R_1$ and $R_2$, then
they will not commute for almost any choice of $R_1$ and $R_2$.  

We make the choices 
\[
R_1 = \left( \begin{array}{ccc} 1 & 0 & 0 \\
0 & {{1}\over{2}} & {{\sqrt{3}}\over{2}} \\
0 & -{{\sqrt{3}}\over{2}} & {{1}\over{2}} \end{array} \right), 
R_2 = \left( \begin{array}{ccc} {{\sqrt{2}}\over{2}} & 0 & {{\sqrt{2}}\over{2}} \\
0 & 1 & 0 \\
-{{\sqrt{2}}\over{2}} & 0 & {{\sqrt{2}}\over{2}} \end{array} \right). \]
This means that $R_3$ and $R_4$ are the following:
\[
R_3 = \left( \begin{array}{ccc} \sqrt{{2}\over{5}} & \sqrt{{3}\over{5}} & 0 \\
-\sqrt{{3}\over{5}} & \sqrt{{2}\over{5}} & 0 \\
0 & 0 & 1 \end{array} \right), 
R_4 = \left( \begin{array}{ccc} 1 & 0 & 0 \\
0 & \sqrt{{5}\over{8}} & -\sqrt{{3}\over{8}} \\
0 & \sqrt{{3}\over{8}} & \sqrt{{5}\over{8}} \end{array} \right). \]

Verification that the approprate pairs of matrices in Table \ref{sl3} do 
not commute is done by explicitly calculating their commutators.  These
are relegated to a separate appendix which is published on the author's web site.

\begin{thm}
For any closed hyperbolic surface group $\Gamma$ and $n\ge 3$, 
there are infinitely many conjugacy classes of discrete, faithful representations 
of $\Gamma$ into $SL(n,\R)$.
\end{thm}

The only thing remaining to prove in this theorem is that we can obtain infinitely
many conjugacy classes.  This is not immediate, since representations of $A(C_5)$ which
are not conjugate might still yield representations of a $\pi_1(\Sigma)$ which are conjugate.
However, the flexibility of our construction proves that you can in fact guarantee 
different conjugacy classes for the surface groups.  

Take any discrete, faithful representation $\sigma: A(C_5) \rightarrow G$ given by 
our construction.  Pick $x \in X = G/K$.  Since $\sigma$ is discrete, we know that
for any $r>0$, the set $B_r = \{a \in A(C_5) | d(\sigma(a)x, x) < r\}$ is finite.  

Let $a_1, a_2, \ldots a_5$ be the generators of $A(C_5)$.  
Then the representation $\sigma_n$ given by $\sigma_n (a_i) = \sigma(a_i)^n$ will be
discrete and faithful.  The image of $\sigma_n$ will be the image of$\sigma$ restricted to the
subgroup $A_n$ of $A(C_5)$ generated by the $a_i^n$;  we can choose $n$ large enough 
so that $A_n \cap B_r = \{ e \}$, since $B_r$ is finite.  

Thus for every $r$, we can 
find a discrete, faithful representation $\tau: A(C_5)\rightarrow G$ such that 
$\min \{ d (x, \tau(a)x) | a \in A(C_5) - \{e\} \} > r$.  

Now suppose that we have a finite set of representations $\tau_1 \ldots \tau_n: \pi_1(\Sigma) \rightarrow G$. 
Take some $a \in \pi_1(\Sigma) \subset A(C_5)$ and let $r = \max d(x,\tau_i(a)x)$.  
If $\tau: A(C_5) \rightarrow G$ is a discrete, faithful representation 
such that $\min \{ d (x, \tau(a)x) | a \in A(C_5) - \{e\} \} > r$, then $\tau$ restricted to 
$\pi_1(\Sigma)$ cannot be conjugate to any of the $\sigma_i$.

\subsection{Other Groups}

The example of $SL(3,\R)$ serves as a model for embeddings of Artin groups on cycles into other Lie groups of rank 2.

Crisp and Wiest's result about embedding surface groups holds only for $A(C_5)$, 
so if the cycle is larger than $C_5$ we would not be able to embed all surface groups into $G$ via this Artin group 
technique. 
However, we would still be able to embed certain surface groups, thanks to the 
result of Droms, Servatius and Servatius.

For instance, we have the following:

\begin{thm}\label{so32}
There are infinitely many conjugacy classes of discrete, faithful representations 
of $A(C_6)$ and $\pi_1(\Sigma_{17})$ into $SO(3,2)$.
\end{thm}

Here the symmetric space is $X=SO(3,2)/SO(3)\times SO(2)$, again of rank 2.  

Let $p_0$ be the identity coset $K = SO(3) \times SO(2)$.  We wish to
find a collection of six geodesics through $p_0$ satisfying the conditions of 
Theorem \ref{mainthm}.

Let $\mathfrak{a}$ be the
abelian subalgebra of $\mathfrak{g}$ generated by the two elements
\[
Y_0 = \left( \begin{matrix}
0 & 0 & 0 & 1 & 0 \\
0 & 0 & 0 & 0 & 0 \\
0 & 0 & 0 & 0 & 0 \\
1 & 0 & 0 & 0 & 0 \\
0 & 0 & 0 & 0 & 0 \end{matrix} \right), 
Y_1 = \left( \begin{matrix}
0 & 0 & 0 & 0 & 0 \\
0 & 0 & 0 & 0 & 1 \\
0 & 0 & 0 & 0 & 0 \\
0 & 0 & 0 & 0 & 0 \\
0 & 1 & 0 & 0 & 0 \end{matrix} \right). \]
Then $\exp(\mathfrak{a})p_0$ is a flat $F_{01}$ in $X$, and $\gamma_0(t) = \exp(tY_0)p_0$ and
$\gamma_1(t) = \exp(tY_1)p_0$ are singular geodesics in $F_{01}$ and are perpendicular.  The 
other singular geodesics in $F_{01}$ through $p_0$ are the angle bisectors $\lambda_0$, $\lambda_1$ 
of the right angles formed by $\gamma_0$ and $\gamma_1$.  Thus if $T_0 = \exp(aY_0)$ and $T_1 = \exp(aY_1)$, 
$T_0T_1$ and $T_0T_1\inv$ are axial isometries of $X$ whose axes are these singular geodesics $\lambda_0$ and
$\lambda_1$.  We choose $a = \log (2 + \sqrt{3})$, so 
\[
T_0 = \left( \begin{matrix}
{2} & 0 & 0 & \sqrt{3} & 0 \\
0 & 1 & 0 & 0 & 0 \\
0 & 0 & 1 & 0 & 0 \\
\sqrt{3} & 0 & 0 & 2 & 0 \\
0 & 0 & 0 & 0 & 1
\end{matrix} \right); 
T_1 = \left( \begin{matrix}
1 & 0 & 0 & 0 & 0 \\
0 & {2} & 0 & 0 & \sqrt{3}  \\
0 & 0 & 1 & 0 & 0 \\
0 & 0 & 0 & 1 & 0  \\
0 & \sqrt{3} & 0 & 0 & {2}  
\end{matrix} \right).
\]

\begin{figure}[htbp]
  \begin{center}
    \[
    \begin{xy}
    	0;<3cm,3cm> **@{-}
    	?>*@{>} ?>*!/-4mm/{T_0T_1}
    	,<0cm,1.5cm>;<3cm,1.5cm> **@{-}
    	?>*@{>} ?>*!/-4mm/{T_0}
    	,<0cm,3cm>;<3cm,0cm> **@{-}
    	?>*@{>} ?>*!/-4mm/{T_0T_1^{-1}}
    	,<1.5cm,0cm>;<1.5cm,3cm> **@{-}
    	?>*@{>} ?>*!/-4mm/{T_1}
    \end{xy}
\]
   \caption{Singular geodesics in a flat in $SO(3,2)/SO(3)\times SO(2)$}
  \end{center}
\end{figure}
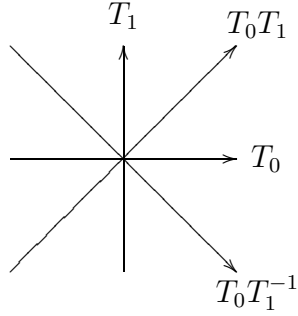

Let $A_0$ and $A_1$ be the subgroups of $K$ fixing $\gamma_0$ and $\gamma_1$, respectively.
There are isomorphisms $\tau_i$ between $SO(2)$ and each of the $A_i$; they are defined by 
\[
\tau_0: \left(\begin{matrix} a & b \\ -b & a \end{matrix} \right) \in SO(2) \mapsto \left( \begin{matrix}
1 & 0 & 0 & 0 & 0 \\
0 & a & b & 0 & 0 \\
0 & -b & a & 0 & 0 \\
0 & 0 & 0 & 1 & 0 \\
0 & 0 & 0 & 0 & 1 \end{matrix} \right) \in A_0 \]\[
\tau_1: \left(\begin{matrix} a & b \\ -b & a \end{matrix} \right) \in SO(2) \mapsto \left( \begin{matrix}
a & 0 & b & 0 & 0 \\
0 & 1 & 0 & 0 & 0 \\
-b & 0 & a & 0 & 0 \\
0 & 0 & 0 & 1 & 0 \\
0 & 0 & 0 & 0 & 1 \end{matrix} \right) \in A_1 \]

Copying our proof of Theorem \ref{sl3}, if we can find $R_1, R_3, R_5 \in A_1$ and $R_2, R_4 \in A_0$ such that 
$R_5R_4R_3R_2R_1 \in A_0$, then $\gamma_0$, $\gamma_1$, $\gamma_2 = R_5\gamma_0$,
$\gamma_3 = R_5R_4\gamma_1$, $\gamma_4 = R_5R_4R_3\gamma_0$, $\gamma_5 = R_5R_4R_3R_2\gamma_1$ are six 
geodesics which are adjacent in a 6-cycle pattern.  

Since $\tau_0$ and $\tau_1$ have image only in the first factor of $K = SO(3) \times SO(2)$,
the condition that $R_5R_4R_3R_2R_1 \in A_0$ is essentially just a requirement that 
a product of rotations in $\R^3$ be a rotation around the $x$-axis.  Just like in our discussion in the
case of $SL(3,\R)$, we see that given any choice of $R_1, R_2, R_3$, we 
know that we can find a rotation $R_4$ around the $x$-axis and a rotation $R_5$ around the $y$-axis such that 
the product $R_5R_4R_3R_2R_1$ fixes the $x$-axis.

As was the case in Theorem \ref{sl3}, we only need to check non-adjacency for one particular 
set of choices of $R_i$.  We choose 
\[
R_1 = \tau_1 \left( \begin{matrix} 
{{3}\over{5}} & {{4}\over{5}} \\
-{{4}\over{5}} & {{3}\over{5}} \end{matrix} \right), 
R_2 = \tau_0 \left( \begin{matrix} {{4}\over{5}} &  {{3}\over{5}} \\
-{{3}\over{5}} & {{4}\over{5}} \end{matrix} \right),
R_3 = \tau_1 \left( \begin{matrix} 
-{{31}\over{481}} & -{{480}\over{481}} \\
{{480}\over{481}} & -{{31}\over{481}} \end{matrix} \right). \]

This forces $R_4$ and $R_5$ to be the following:
\[
R_4 = \tau_0 \left( \begin{matrix}  {{4}\over{5}}   & {{3}\over{5}}  \\
-{{3}\over{5}}  & {{4}\over{5}} 
\end{matrix} \right), 
R_5 = \tau_1 \left( \begin{matrix} {{3}\over{5}} & {{4}\over{5}} \\
-{{4}\over{5}} &  {{3}\over{5}} \end{matrix} \right). 
\]

Axial isometries with the $\gamma_i$ as axes are $T_0$, $T_1$, $R_5T_0R_5\inv$, $R_5R_4T_1 R_4\inv R_5\inv$, 
$R_5R_4R_3T_0R_3\inv R_4\inv R_5\inv$, $R_5R_4R_3R_2T_1R_2\inv R_3\inv R_4\inv R_5\inv$.  In addition,
the extra singular geodesics in the six flats will have axial isometries listed in the
following table, where  $F_{ij}$ denotes the flat containing $\gamma_i$ and $\gamma_j$.   

\begin{table}[h]\label{so32table}
\begin{center}
\begin{tabular}{c|c}
	Flat & Isometries \\
        \hline
        $F_{01}$ & $T_0T_1$ \\
         \  & $T_0T_1\inv$ \\
	\hline
	$F_{12}$ & $R_5T_0T_1R_5\inv$ \\
         \  & $R_5T_0T_1\inv R_5\inv$ \\
       \hline
       	$F_{23}$ & $R_5R_4T_0T_1R_4\inv R_5\inv$ \\
         \  & $R_5R_4T_0T_1\inv R_4\inv R_5\inv$ \\
         \hline
	 	$F_{34}$ & $R_5R_4R_3T_0T_1R_3\inv R_4\inv R_5\inv$ \\
         \  & $R_5R_4R_3T_0T_1\inv R_3\inv R_4\inv R_5\inv$ \\
         \hline
	 	$F_{45}$ & $R_5R_4R_3R_2T_0T_1R_2\inv R_3\inv R_4\inv R_5\inv$ \\
         \  & $R_5R_4R_3R_2T_0T_1\inv R_2\inv R_3\inv R_4\inv R_5\inv$ \\
         \hline
	 	$F_{50}$ & $R_5R_4R_3R_2R_1T_0T_1R_1\inv R_2\inv R_3\inv R_4\inv R_5\inv$ \\
         \  & $R_5R_4R_3R_2R_1T_0T_1\inv R_1\inv R_2\inv R_3\inv R_4\inv R_5\inv$ \\
\end{tabular}
\caption{Axial isometries by flat}
\end{center}
\end{table}

Again, checking the non-adjacency requirements (of which there are 93) is an explicit computation of commutators of these matrices, and is published in a separate appendix on the author's web site.  

The proof that there are infinitely many conjugacy classes proceeds in exactly the same way as
in Theorem \ref{sl3}.

We have been unable to find a configuration of five maximal flats in $X$ 
that satisfy all of the requirements of Theorem \ref{mainthm}.  It is possible to find flats satisfying the first two conditions, but in every configuration analyzed, there are always other geodesics in the flats which are adjacent.  Changing the Artin group to $A(C_6)$ allows for more freedom in choosing the rotations around each singular geodesic and allows us to avoid this problem. 

It seems likely that Artin groups on cycles of 6 or more should be flexible enough to allow representations into most semisimple Lie 
groups, barring a few obstructions.  

For instance, suppose there are only two singular geodesics in a flat through a given point.  This is the case when
$X = SU(2,n)/S(U(2) \times U(n))$, for instance.  In this situation, 
it will be impossible to embed $A(C_n)$ into the associated Lie group
via this method if $n$ is odd, since there will have to be two adjacent
geodesics which are conjugates of the same geodesic in the original flat.

Another obstruction occurs for non-irreducible symmetric spaces. If $X$ is  a product of rank one symmetric spaces $X_1 \times X_2$, no
$A(C_n)$ (for $n>4$) will embed faithfully by this method.  In this case, 
if $p_0 = (p_1, p_2)$, singular geodesics through $p_0$ are either $\{\eta_1\} \times \{p_2\}$
or $\{p_1\} \times \{\eta_2\}$, where $\eta_i$ is a geodesic through $p_i$ in $X_i$, and any
geodesic of the first type will be adjacent to any geodesic of the second, so there is no
hope of satisfying the conditions of Theorem \ref{mainthm} for graphs which are not
complete and bipartite.  

The general conjecture looks like this:  

\begin{conj}
If $X = G/K$ is an irreducible symmetric space of rank at least two, and given any maximal flat $F$
and point $p_0 \in F$, there are at least three singular geodesics in $F$ passing through $p_0$, then
there are infinitely many conjugacy classes of discrete, faithful representations of $A(C_n)$ into $G$ for
any $n\ge 6$.

If $X = G/K$ is an irreducible symmetric space of rank two, and given any maximal flat $F$
and point $p_0 \in F$, there are exactly two singular geodesics in $F$ passing through $p_0$, then 
there are infinitely many conjugacy classes of discrete, faithful representations of $A(C_{2n})$ into $G$,
for any $n \ge 3$.  
\end{conj}

It would be helpful and interesting to find general methods that
would allow for attacking this problem in ways other than a 
case-by-case analysis.

\section{Lattices}

We can also embed Artin groups inside some lattices in higher rank Lie groups.
Recall that a {\emph{lattice}} in $G$ is a discrete subgroup $\Gamma\subset G$ such that
$\Gamma \backslash G$ has finite volume.  Embeddings of 
surface groups into lattices give examples of essential surfaces inside locally
symmetric spaces.  

\begin{thm}\label{sl5z}
There are infinitely many conjugacy classes of representations of $A(C_5)$ into 
$\Gamma = SL(5, \Z)$.  
\end{thm}

We will find five geodesics to apply the main theorem to.  

Let $p_0$ be the identity coset  in the symmetric space $X=SL(5,\R)/SO(5,\R)$.
Fix an integer $n\ge 2$, and let $A_i$ be the following elements of $\Gamma$:
\[
A_1 = \left( \begin{matrix} 
n & n-1 & 0  & 0 & 0 \\
n+1 & n & 0 & 0 & 0 \\
0 & 0 & 1 & 0 & 0 \\
0 & 0 & 0 & 1 & 0 \\
0 & 0 & 0 & 0 & 1 \end{matrix} \right),
A_2 = \left( \begin{matrix} 
1 & 0 & 0 & 0 & 0 \\
0 & 1 & 0 & 0 & 0 \\
0 & 0 & n & n-1 & 0 \\
0 & 0 & n+1 & n & 0 \\
0 & 0 & 0 & 0 & 1 \end{matrix} \right), \] \[
A_3 = \left( \begin{matrix} 
n & 0 & 0 & 0 & n+1 \\
0 & 1 & 0 & 0 & 0 \\
0 & 0 & 1 & 0 & 0 \\
0 & 0 & 0 & 1 & 0 \\
n-1 & 0 & 0 & 0 & n \end{matrix} \right),
A_4 = \left( \begin{matrix} 
1 & 0 & 0 & 0 & 0 \\
0 & n & n-1 & 0  & 0 \\
0 & n+1 & n & 0 & 0 \\
0 & 0 & 0 & 1 & 0 \\
0 & 0 & 0 & 0 & 1 \end{matrix} \right), \] \[
A_5 = \left( \begin{matrix} 
1 & 0 & 0 & 0 & 0 \\
0 & 1 & 0 & 0 & 0 \\
0 & 0 & 1 & 0 & 0 \\
0 & 0 & 0 & n & n-1 \\
0 & 0 & 0 & n+1 & n \end{matrix} \right).
\]

As isometries of $X$, these are all axial hyperbolic isometries.  
We see easily that $A_i$ and $A_j$ commute if and only if $i-j \equiv \pm 1$ (mod 5), 
so the geodesics $\gamma_i (t) = \exp(t \log A_i)SO(5,\R)$ mimic the graph $C_5$ in the appropriate way.  

The symmetric space is of rank 4, so the requirement that the set of singular directions in the span of 
$\log A_i$ and $\log A_j$ be finite is non-trivial.  It is relatively easy to check, however.  The maximal
flat containing both $\gamma_1$ and $\gamma_2$ can be taken to the maximal flat $F$ of diagonal matrices via 
some element $k_{12}$ of $K$ which conjugates $A_1$ and $A_2$ to the matrices
\[
B_1 = \left( \begin{matrix} 
\lambda & 0 & 0  & 0 & 0 \\
0 & \lambda\inv & 0 & 0 & 0 \\
0 & 0 & 1 & 0 & 0 \\
0 & 0 & 0 & 1 & 0 \\
0 & 0 & 0 & 0 & 1 \end{matrix} \right), 
B_2 = \left( \begin{matrix} 
1 & 0 & 0  & 0 & 0 \\
0 & 1 & 0 & 0 & 0 \\
0 & 0 & \lambda & 0 & 0 \\
0 & 0 & 0 & \lambda\inv & 0 \\
0 & 0 & 0 & 0 & 1 \end{matrix} \right), \]
where $\lambda$ and $\lambda\inv$ are the eigenvalues of the matrix $\left( \begin{matrix} 
n & n-1 \\
n+1 & n \end{matrix} \right)$.  
In the Lie algebra $\mathfrak{a}$ of diagonal traceless matrices, the span corresponds to the subalgebra $\mathfrak{b}$
of all elements of the form $$\left( \begin{array}{ccccc} 
a & 0 & 0  & 0 & 0 \\
0 & -a & 0 & 0 & 0 \\
0 & 0 & b & 0 & 0 \\
0 & 0 & 0 & -b & 0 \\
0 & 0 & 0 & 0 & 0 \end{array} \right).$$
Since the only singular directions in $\mathfrak{a}$ are those with equal entries on the diagonal,
we see that there are only 4 singular directions in $\mathfrak{b}$, namely those where $a=0, b=0, 
a=b,$ or $a=-b$.  

It remains to check the non-adjacency requirements.  Let $F_{ij}$ denote the two-plane containing 
$\gamma_i$ and $\gamma_j$.  Then the singular geodesics through $p_0$ in $F_{ij}$ are precisely the axes of the
isometries $A_i$, $A_j$, $A_iA_j$, and $A_iA_j\inv$.  There are 125 pairs of geodesics that are
supposed to be non-adjacent, but by symmetry arguments we can focus on those involving one geodesic
from $F_{12}$ and one geodesic from either $F_{23}$ or $F_{34}$.  

Since $A_1$ commutes with $A_2$ but not with $A_3$, it will not commute 
with either $A_2A_3$ or $A_2A_3\inv$.  Neither $A_1A_2$ nor $A_1A_2\inv$ 
commute with $A_3$, so they also will not commute with $A_2A_3$ and $A_2A_3\inv$.  
Thus all of the pairs of geodesics in $F_{12}$ and $F_{23}$ behave the way they ought.

Checking that no geodesic in $F_{12}$ is adjacent to one in $F_{34}$ (save $\gamma_2$ being adjacent $\gamma_3$) 
requires a few more explicit matrix calculations, which we relegate to a separate appendix on the author's web site.  

Thus the geodesics satisfy the requirements of Theorem \ref{mainthm}, and so high enough
powers of the $A_i$ will give a discrete and faithful subgroup of $SL(5,\Z)$ isomorphic to $A(C_5)$.  

Note that we could get the image of the representation to lie in any congruence subgroup that we desire,
since for any prime $p$, arbitrarily high powers of the $A_i$ will be congruent to the identity modulo $p$.

\addcontentsline{toc}{chapter}{References}

%
%

\begin{flushleft}
\sc{Department of Mathematics \\
Haverford College \\
Haverford, PA  19104} \\
\tt{swang@haverford.edu}
\end{flushleft}

\end{document}